\documentclass[12pt,twoside]{article}

\usepackage{amsmath,amsfonts,amssymb}
\usepackage{graphicx}

\newcommand{\law}[1]{\mathcal{L}\bigl(#1\bigr)}

\newcommand{\eqd}{\stackrel{d}{=}}
\newcommand{\pr}{{\bf P}}

\newtheorem{theo}{\sc Theorem}

\begin{document}

\title{On the distribution of the Brownian motion process on its way to hitting zero}

\author{\sc Konstantin Borovkov$^1$\\
\small {\em Dept. of Maths \&\ Stats, University of Melbourne,}\\
\small {\em Parkville 3010, Australia}\\
\small email: {\tt borovkov@unimelb.edu.au}. }

\date{}

\maketitle

\footnotetext[1]{Research supported by ARC Discovery Grant DP0880693.}

\begin{abstract}
We present functional versions of recent results on the univariate distributions of the
process $V_{x,u} = x + W_{u\tau(x)},$ $0\le u\le 1$, where $W_\bullet$ is the standard
Brownian motion process,  $x>0$ and $\tau (x) =\inf\{t>0 :\, W_{t}=-x\}$.
\end{abstract}

{\noindent\small AMS 2000 Subject Classification:   60J65\\
 \vspace{5mm}%
 }

Let $\{W_t\}_{t\ge 0}$ be the standard univariate Brownian motion process and, for
$x>0$,
\[
W_{x,t}:= x+W_t, \quad t\ge 0,\qquad \tau (x) :=\inf\{t>0 :\, W_{x,t}=0\}.
\]
As is well known, $\tau (x)$ is a proper random variable with density
\begin{equation}
\label{PDF}
p_x(t)=\frac{x e^{-x^2/2t}}{\sqrt{2\pi t^3}} , \qquad t>0,
\end{equation}
so one can introduce
\[
V_{x,u} := W_{x,u\tau(x)}, \qquad 0\le u\le 1.
\]
These random variables were studied in the recent paper~\cite{KlCh}, where it was shown
(Theorem~1.1) that, for any fixed $u\in (0,1),$   $V_{x,u}$  has density
\begin{align}
p_{x,u} (y) & : = \frac{d}{dy} \pr (V_{x,u}\le y)   &\notag\\
\label{Fima1}
   & = \frac{4 \sqrt{u(1-u)}\, x y^2}{\pi  \left(u y^2+(1-u)
   (y-x)^2\right) \left(u y^2+(1-u) (y+x)^2\right)},
   &  y>0,  \\
\label{Fima1a}
   &\sim \frac{4 x\sqrt{u(1-u)}}{\pi y^2} &   \mbox{as $y\to\infty$ }
\end{align}
(here and in what follows, $a\sim b$ means that $a/b\to 1$). Representation
\eqref{Fima1} implies, in particular, that, for any fixed $u\in [0,1]$, one has
\begin{equation}
\label{Fima2}
  V_{x,u} \eqd x V_{1,u}.
\end{equation}
Using a direct tedious calculation, it was also demonstrated in Section~3 of~\cite{KlCh}
that, for a  fixed $u\in (0,1)$, the density $p_{x,u}$ coincides with that of a ``scaled
Brownian excursion at the corresponding time, averaged over its length". The
mathematical formulation of that result was given by formula~(3.3) in~\cite{KlCh}  that
can be rewritten as follows. Let $\{R^T_{x,t}\}_{t\le T}$ be a three-dimensional Bessel
bridge of length $T$ pinned at $x$ at time $t=0$ and at $0$ at time~$t=T$, which is
independent of our process $W_{x,\bullet}$ (and hence of $\tau(x)$). Recall that one can
represent the process as
\begin{equation}
\label{My1} R^T_{x,t} = \bigl\|   W^{(3)}_{x,t}  - tT^{-1} W^{(3)}_{x,T} \bigr \|,
\qquad 0\le t \le T,
\end{equation}
where
\begin{equation}
\label{My2}
W^{(3)}_{x,t} = (x,0,0) +  W^{(3)}_{t} , \qquad t\ge 0,
\end{equation}
$W^{(3)}_{\bullet}$ being   a standard three-dimensional Brownian motion process and
$\|\cdot \|$   the Euclidean norm in~${\mathbb R}^3$. The above-mentioned formula
from~\cite{KlCh} is equivalent to the assertion that, for any fixed $u\in [0,1],$ one
has
\begin{equation}
\label{Fima3}
V_{x,u} \eqd R^{\tau (x)}_{x,u \tau (x)}.
\end{equation}
Note that $R^T_{x,\bullet}$ is not exactly an excursion (an excursion returns to the
same point where it started) but, rather, a time-reversed Brownian meander (see e.g.\
p.63 in~\cite{BoSa}), and that on the right-hand side of~\eqref{Fima3} it is averaged
not over its length, but rather of that of an independent version of $W_{x,\bullet}$ on
its way to hitting zero.

Observe also that, due to the self-similarity of the Brownian motion process,
representation \eqref{My1}--\eqref{My2} implies that
\[
R^T_{x,\bullet T}  \eqd T^{1/2} R_{T^{-1/2} x,\bullet}
\]
(as processes in $C[0,1]$), where we put $R_{x,t}:=R^1_{x,t}.$

The main aim of the present note is to give  simple proofs to functional versions
of~\eqref{Fima2} and~\eqref{Fima3} (that had ``remained elusive", as was noted
in~\cite{KlCh}).

\begin{theo}
\label{Th1}
For any $x>0$,
\begin{equation}
\label{My3_0} \{V_{x,u}\}_{u\le 1} \eqd  \{x V_{1,u}\}_{u\le 1}.
\end{equation}
Furthermore, there exists a regular version of the conditional distribution of
$V_{1,\bullet}$ in $C[0,1]$ given $\tau (1)=T$ that coincides with the law of $T^{1/2}
R_{T^{-1/2},\bullet}$, and therefore, if  $\tau\eqd \tau (1)$ is independent of the
Brownian motion process from representation \eqref{My1}--\eqref{My2}, then one has
\begin{equation}
\label{My3}
\{V_{1,u}\}_{u\le 1} \eqd  \{  \tau^{1/2} R_{\tau^{-1/2},u}\}_{u\le 1}.
\end{equation}
\end{theo}

Before proceeding to prove Theorem~\ref{Th1}, we will make the following observation.
That the asymptotic behaviour \eqref{Fima1a} holds is most natural in view of
representation~\eqref{My3} and the form \eqref{PDF} of the density of $\tau(x)$, the
latter implying that $\tau^{1/2}$ has the density
\[
\sqrt{\frac2\pi}  y^{-2} e^{-1/2y^2}, \qquad y>0.
\]
Moreover, since $W_{x,\bullet}$ is a martingale with $W_{x,0}=x$, from the optional
sampling theorem one immediately obtains that, using notation
\[
\hat X_{t} := \sup_{s\le t} X_s, \qquad \check{X_t}  := \inf_{s\le t} X_s
\]
for the maximum and minimum processes of  $X_\bullet$, one has
\begin{equation}
\label{My3a}
\pr (\hat V_{x,1} > y) =\pr \bigl(  \hat W_{x,\tau (x)} > y\bigr) =x y^{-1}, \qquad y>x,
\end{equation}
with the same asymptotic behaviour $ \varpropto y^{-2}$ as
$y\to \infty$ for the density.

Finally, denoting by $q_{x,y,u} (z)$ the density of a random variable following the
conditional distribution $\law{y^{-1} V_{x,u}\, | \, \hat V_{x,1} > y}$ of $y^{-1}
V_{x,u}$ given $\hat V_{x,1} > y$ (here and in what follows, $\law{  X\, | \, C}$
denotes the conditional distribution of the random element $X$ in the respective
measureable space given condition~$C$), it is easily seen from \eqref{Fima1},
\eqref{Fima1a} and \eqref{My3a} that, for any fixed $x,z>0$ and $u\in (0,1)$,
\[
q_{x,y,u} (z) \sim \frac{4 \sqrt{u(1-u)}}{\pi z^2} \qquad \mbox{as}\quad  y\to\infty.
\]

\noindent {\sc Proof of Theorem~\ref{Th1}.} First we observe that
\begin{equation}
\label{tilda}
W_{x,t} = x(1+ x^{-1}W_t) = x \widetilde W_{1,tx^{-2}}, \qquad t\ge 0,
\end{equation}
where $\widetilde W_{1,\bullet}$ is a Brownian motion process starting at~1. All
quantities related to this process we will label with tilde. As $\tau (x)$ is the first
time the LHS of \eqref{tilda} turns into zero, we see that $\widetilde \tau (1) =\tau
(x) x^{-2}.$ Therefore
\[
V_{x,u} = W_{x, u\tau (x)} = x \widetilde W_{1,u \widetilde \tau (1)} = x \widetilde
V_{1,u}, \qquad u\in [0,1],
\]
which proves~\eqref{My3_0}. So from now on, we can assume without loss of generality
that~$x=1$.

Next let, for some functions $f_j  \in C[0,1]$ and numbers $r_j>0,$ $j=1,2,\ldots, n,$
\[
A:=\bigcap_{j\le n} \{ f\in C[0,1]: \, \|f - f_j\|<r_j\}
\]
be a finite intersection of open balls in $C[0,1]$ ($\|\cdot\|$ stands for the uniform
norm). For $T,h,\delta >0$, put
\[
A_T:= \{f (\mbox{\small$\bullet$} /T):\, f\in A\}\subset C[0,T], \qquad
 \varepsilon ( \delta) := \max_{j\le n} \omega_{f_j} (\delta),
\]
where $\omega_f (\delta):=\sup_{0\le s <t \le s+\delta \le 1} |f(s) - f(t)|$ is the
continuity modulus of the function~$f$. Finally, we denote by $A_{T}^{\varepsilon
(h/T)}$   the $ \varepsilon (h/T)$-neighbourhood of $A_T$ (in the uniform topology on~$
C[0,T]$) and introduce the event
\[
B_{T,h} := \bigl\{ \{W_{1,t}\}_{t\in [0,T]}  \in A_{T}^{\varepsilon (h/T)}\bigr\}.
\]

Now, using the Markov property and the well-know relations
\[
\pr \bigl(  \check W_{T+h} <0 \,|\, W_T=  y\bigr) =  2\overline \Phi (yh^{-1/2}),
  \quad
\pr  \bigl(\check W_{1,T} >0  \,|\,   W_{1,T}=  y\bigr) = 1 - e^{-2y/T}, \quad y>0,
\]
where $\overline\Phi=1-\Phi,$ $\Phi $ being the standard normal distribution function,
we have
\begin{align}
\pr & \bigl(V_{1,\bullet}  \in A ,\, \tau (1) \in (T, T+h)\bigr)
 \le \pr \bigl( B_{T,h},\, \tau (1) \in (T, T+h)\bigr)
  \notag\\
 & = \int_0^\infty \pr \bigl(B_{T,h} ,\, \tau (1) \in (T, T+h) \,|\, W_{1,T}  =
 y\bigr)  \pr( W_{1,T}\in dy)
 \notag \\
 & = \int_0^\infty \pr \bigl(B_{T,h}, \check W_{1,T} >0,  \check W_{1,T+h} <0 \,|\, W_{1,T}=
 y\bigr)  \pr( W_{1,T}\in dy)
 \notag\\
 & = \int_0^\infty \pr \bigl(B_{T,h} , \check W_{1,T} >0   \,|\, W_{1,T}=
 y\bigr) \pr \Bigl(  \min_{t\in [T, T+h]} W_{1,t}  <0 \,\Big|\, W_{1,T}=
 y\Bigr)  \pr( W_{1,T}\in dy)
 \notag\\
& = \int_0^\infty \pr \bigl(B_{T,h} \,|\, \check W_{1,T} >0 ,   W_{1,T}=
 y\bigr) \pr  \bigl(\check W_{1,T} >0  \,|\,   W_{1,T}=
 y\bigr)   2\overline \Phi (yh^{-1/2})  \pr( W_{1,T}\in dy)
 \notag\\
& = 2 \int_0^\infty \pr \bigl(B_{T,h} \,|\, \check W_{1,T}  >0 ,   W_{1,T}=
 y\bigr)   \bigl(1 - e^{-2y/T}\bigr) \overline \Phi (yh^{-1/2})  \pr( W_{1,T}\in dy)
 \notag\\
 &  = (4+ o(1))h^{1/2}  \int_0^{h^{1/4}} \pr \bigl(B_{T,h} \,|\, \check W_{1,T}  >0 ,   W_{1,T}=
 y\bigr) g_T (yh^{-1/2}) \, dy + o(h)
 \label{dlina}
\end{align}
as $h\downarrow 0,$ where
\[
g_T (u) = \frac1{\sqrt{ 2\pi}}  \,  u T^{-3/2} e^{-1/(2T)} \overline \Phi (u),
 \qquad u>0,
\]
and we used Mills ratio to infer that $\int_{h^{1/4}}^\infty = o(h)$

Next we will show that the   probability in the last integrand in~\eqref{dlina}
converges to the respective probability for the Brownian meander process as
$y\downarrow 0$.

Recall that the Brownian meander process $\{W^\oplus_s\}_{s\le 1}$ can be defined as
follows (see e.g.~\cite{DuIgMi} or p.64 in~\cite{BoSa}): letting $\zeta :=\sup\{t\le
1:\, W_t =0\}$ be the last zero of the Brownian motion in $[0,1]$, we set
\[
W^\oplus_s := (1 - \zeta)^{-1/2} \big|W_{\zeta + (1-\zeta) s}\big|, \qquad 0\le s\le 1.
\]
This is a continuous nonhomogeneous Markov process whose transition density can be found
e.g.\ in~\cite{DuIgMi} (relations  (1.1) and~(1.2)). It is known that the conditional
version of the process pinned at $x>0$  at time $s=1$ coincides in distribution with the
three-dimensional Bessel process starting at zero and also pinned at $x$ at time $s=1$
(see e.g.\ p.64 in~\cite{BoSa}), which can be written as
\[
\law{ \{W^\oplus_s\}_{s\le 1} \,| \, W^\oplus_1 = x}
 = \law{ \{ \|W^{(3)}_s\|\}_{s\le 1} \,| \,  \|W^{(3)}_1\| =x}.
\]
It is not hard to deduce from here, the spherical symmetry of the Brownian motion
process $W^{(3)}_\bullet$ and representation \eqref{My1}--\eqref{My2} above that
\begin{equation}
\label{My3b}
\law{ \{W^\oplus_s\}_{s\le 1} \,| \, W^\oplus_1 = x}
 = \law{ \{ \|  W^{(3)}_{x,1-s} - (1-s) W^{(3)}_{x,1} \|\}_{s\le 1}}
 = \law{ \{    R_{x,1-s} \}_{s\le 1}} .
\end{equation}

An alternative insightful  interpretation of the Brownian meander is given by the fact
that its distribution (in $C[0,1]$) coincides with the weak limit of conditional
distributions of $W_\bullet$ conditioned to stay above $-\varepsilon \uparrow 0:$
\[
\law{ \{W^\oplus_s\}_{s\le 1}} = \mbox{w-}\!
 \lim_{\varepsilon\downarrow 0}
  \law{\{W_{s}\}_{s\le 1} \,| \, \check W_{1} >- \varepsilon}
\]
(Theorem~(2.1) in~\cite{DuIgMi}; $\mbox{w-}\! \lim $ stands for the limit  in weak
topology). A conditional version of a result of this type is used in the calculation
displayed in~\eqref{dlina1} below.

Now return   to the   probability in the integrand in the last line in~\eqref{dlina} and
recall the well-known property of Brownian bridges that conditioning a Brownian motion
on its arrival at a point~$y\neq 0$ at  time~$T$ is equivalent to conditioning on its
arrival to zero at that time and then adding the deterministic linear trend component
$yt/T$. This implies that, for any $\varepsilon \ge \varepsilon(h/T),$
\begin{align}
\pr &\bigl(B_{T,h} \,|\, \check W_{1,T}  >0 ,   W_{1,T}=
 y\bigr) \notag\\
 & =\pr \bigl( \{W_{1,t} + ytT^{-1}\}_{t\le T} \in A_T^{\varepsilon(h/T)} \,|\, W_{1,T}=0; \, W_{1,s}
> - ysT^{-1}, s\in [0,T]\bigr)\notag\\
 &=
 \pr  \bigl( \{W_{T-t} + ytT^{-1}\}_{t\le T} \in A_T^{\varepsilon(h/T)} \,|\, W_{T}=1; \,   W_{s} > - y(T-s)T^{-1}, s\in
[0,T]\bigr)\notag\\
 &=
\pr  \bigl(\{T^{1/2}W_{1-v} + yv\}_{v\le 1}\in A^{\varepsilon(h/T)} \,|\, W_{1}=
T^{-1/2}; \,
   W_{v} > - yT^{-1/2}(1-v) , v\in [0,1]\bigr)\notag\\
 &\le \pr  \bigl(\{T^{1/2}W_{1-v} + yv\}_{v\le 1}\in A^{\varepsilon } \,|\,
W_{1}= T^{-1/2}; \,
   W_{v} > - yT^{-1/2}(1-v) , v\in [0,1]\bigr)\notag\\
      & \to
    \pr  \bigl(\{T^{1/2}W^\oplus_{1-v}  \}_{v\le 1} \in A^{\varepsilon}\,|\, W^\oplus_{1}= T^{-1/2}
    \bigr)\notag\\
 &     =     \pr  \bigl(\{T^{1/2}R_{T^{-1/2}, s}  \}_{s\le 1} \in A^{\varepsilon} \bigr)
 \label{dlina1}
\end{align}
as $y\downarrow 0,$ where the second last relation follows from the weak convergence
established in Theorem~6 in~\cite{BoDo} (as it is obvious that $A^{\varepsilon}$ has
null boundary w.r.t.\ the limiting distribution) and the last one  follows
from~\eqref{My3b}.

Since  $\varepsilon (h/T)\to 0$ as $h\downarrow 0$, and $A$ has a null boundary under
$\law{\{ T^{1/2} R_{T^{-1/2}, s}  \}_{s\le 1}}$,  we conclude from~\eqref{dlina}
(changing there the variables: $u=yh^{-1/2}$) that
\begin{align}
\limsup_{h\downarrow 0}
  &\,
 \frac{1}{h}\,\pr    \bigl(V_{1,\bullet}  \in A ,\, \tau (1) \in (T, T+h)\bigr)\notag\\
 & \le  \limsup_{h\downarrow 0} 4  \pr  \bigl( T^{ 1/2}R_{T^{-1/2}, \bullet }  \in A \bigr)
   \int_0^{h^{-1/4}}   g_T (u) \, du
 \notag\\
  & =    \pr  \bigl( T^{ 1/2}R_{T^{-1/2}, \bullet }  \in A \bigr) p_1 (T),
 \label{ha}
\end{align}
owing to $ \int_0^\infty   u\overline \Phi (u) \, du = \frac14$ and~\eqref{PDF}.

As the same argument as employed in  \eqref{dlina1} and~\eqref{ha} will also work for
the complement of~$A$, we obtain that
\[
\pr    \bigl(V_{1,\bullet}  \in A ,\, \tau (1) \in (T, T+h)\bigr)
 \sim  \pr  \bigl( T^{ 1/2}R_{T^{-1/2}, \bullet }  \in A \bigr) p_1 (T) h
 \qquad\mbox{as \quad}
 h\downarrow 0.
\]
This relation implies that, for any fixed $0<T_1 <T_2 <\infty,$
\[
\pr   \bigl(V_{1,\bullet}  \in A ,\, \tau (1) \in (T_1, T_2)\bigr)
 =
 \int_{T_1}^{T_2} \pr  \bigl( T^{ 1/2}R_{T^{-1/2}, \bullet }  \in A \bigr)
p_1 (T) \, dT.
\]
Since intersections of finite collections of open balls form  determining classes in
separable spaces (see e.g. Section I.2 in~\cite{Bi}), this completes the proof of the
theorem.

\bigskip\noindent {\bf Acknowledgment.} The author is grateful to the referee whose valuable comments
helped to improve the exposition of the paper.

\end{document}